\documentclass[12pt,twoside]{amsart}
 \usepackage{amssymb}
 \usepackage[all]{xy}
 \usepackage[colorlinks,urlcolor=blue,linkcolor=blue]{hyperref}

 \numberwithin{equation}{section} 
 \newtheorem{theorem}{Theorem}[section]
 
 \newtheorem{lemma}[theorem]{Lemma}
 \newtheorem{corollary}[theorem]{Corollary}
 \theoremstyle{remark}
 \newtheorem{remark}[theorem]{Remark}

 \newcommand{\ra}{\rightarrow}
 \newcommand{\D}{\Delta}
 \newcommand{\wt}{\widetilde}

\begin{document}

 \title{Nonvanishing Betti numbers of weakly chordal graphs}
 \author{Jos\'e Mart\'{\i}nez-Bernal\ \ and\ \ Oscar A. Piza-Morales}
 \address{Departamento de Matem\'aticas, Cinvestav-IPN, M\'exico\\
 \scriptsize A.P. 14-740, Ciudad de M\'exico 07360}
 \email{\{jmb, oapiza\}@math.cinvestav.mx}
 \keywords{Betti numbers, Hochster's formula, Stanley-Reisner ring, weakly chordal graph.}
 \subjclass[2000]{13D02, 13F55}
 \maketitle

\begin{abstract}We show that Kimura's necessary and sufficient condition for the
non-vanishingness of multigraded Betti numbers of chordal graphs can be extended
to the wider class of weakly chordal graphs.
\end{abstract}

\section{Introduction}

Let $R=k[x_1,\ldots,x_n]$ be a polynomial ring with its canonical
$\mathbb{Z}^n$-grading over a field $k$. For a graded finite $R$-module $M$ its
{\it multigraded Betti numbers} are defined as
$\beta_{i\sigma}(M):=\dim_k\text{Tor}_i(M,k)_\sigma$ for all
$\sigma\in\mathbb{Z}^n$. A central problem in commutative algebra is the
understanding of these numbers. Even in the case when $M$ is a quadratic
squarefree monomial ideal of $R$ this is a wide open problem. Here we are
interested in the non-vani\-shin\-gness of Betti numbers for this class of
monomial ideals. Since these ideals are in a natural one-to-one correspondence
with graphs, it can be expected to obtain a non-vanishingness criterion in terms
of the combinatorics of the graph. In~\cite{Zheng} Zheng succeed in to obtain
such a criterion when the graph is a tree. In~\cite{Katzman} Katzman found a
sufficient condition for non-vanishingness in arbitrary graphs. This sufficient
condition was generalized by Kimura in~\cite{Kimura11}, who, moreover, was able
to prove that it is also a necessary one for the class of chordal graphs,
extending in that way Zheng's criterion. In~\cite{Kimura14} Kimura generalized
her sufficient condition given in~\cite{Kimura11}, and our aim here is to show
that her generalized condition is also a necessary one for the wider class of
weakly chordal graphs (Theorem~\ref{mainwchg}). As an application, using a
result of Dalili and Kummini~\cite{DalKumm}, we conclude the well-known fact
that projective dimension and regularity of weakly chordal graphs does not
depend on the characteristic of the base field and we recover the combinatorial
descriptions of these two homological invariants (Corollary~\ref{maincor}).

\section{Preliminaries}

A \emph{simplicial complex} $\D$ on the set of vertices $[n]:=\{1,\ldots,n\}$ is
an inclusionwise closed family of subsets of $[n]$, i.e. $\sigma\in\D$ and
$\tau\subseteq\sigma$ implies $\tau\in\D$; we assume that all the vertices
belongs to $\D$. Elements of $\D$ are called faces. By abuse of notation it is
convenient to identify $\sigma\subseteq[n]$ with the characteristic vector
$\sigma=(\sigma_i)\in\{0,1\}^n$ such that $\sigma_i = 1$ if and only if $i \in
\sigma$; and write $|\sigma|:=\sigma_1+\cdots+\sigma_n$. For
$\sigma\subseteq[n]$ we denote by $\D_\sigma$ the simplicial complex that
results from the restriction of $\D$ to the vertex set $\sigma$.

Given a simplicial complex $\D$, let $I_\D$ denote its \emph{Stanley-Reisner
ideal} and $k[\D]$ its \emph{Stanley-Reisner ring}, i.e. $I_\D=\langle
x_{i_1}\cdots x_{i_r}: \{i_1,\ldots,i_r\}\notin\D\rangle\subset R$ and
$k[\D]=R/I_\D$. Let's also denote by $\wt{H}_i(\D)$ the $i$-th \emph{reduced
homology group} of $\D$ with coefficients in the field $k$. \emph{Hochster's
Formula} is a basic result.

\begin{theorem}\label{MHF} {\rm(Hochster's Formula \cite{H})}
Let $\D$ be a simplicial complex with vertex set $[n]$. For any
$\sigma\subseteq[n]$ we have that $\beta_{i\sigma}(R/I_\D) = \dim
\wt{H}_{|\sigma|-i-1}( \D_\sigma).$\end{theorem}

Let $G=(V,E)$ be a (finite and simple) graph. Two edges $e$ and $f$ of $G$ are
called $3$-\emph{disjoint} if they are vertex disjoint and none endpoint of $e$
is adjacent to an endpoint of $f$. A family $\{B_1,\ldots,B_r\}$ of complete
bipartite subgraphs of $G$ (non-necessarily induced) is called \emph{strongly
disjoint\/} if $V(B_i) \cap V(B_j) = \emptyset$ for all $i \neq j$ and there are
edges $e_1,\ldots, e_r$, with $e_i \in E(B_i)$, which are pairwise $3$-disjoint
in $G$.

A subset $\sigma\subseteq V(G)$ is called an \emph{independent set} if there is
no edge of $G$ with both endpoints in $\sigma$. The \emph{independence complex}
of $G$, denoted $\D(G)$, is the simplicial complex with vertex set the vertices
of $G$ and faces the independent sets of $G$.

Let $I(G)$ denote the {\it edge ideal} of $G$, i.e. $I(G)\subset R$ is the ideal
generated by the monomials $x_ix_j$ such that $\{i,j\}$ is an edge of $G$. A
straightforward verification shows that $I(G)$ is precisely the Stanley-Reisner
ideal of the independence simplicial complex of $G$, i.e. $I(G)=I_{\D(G)}$.
Kimura has given a sufficient condition for the non-vanishingness of Betti
numbers of an arbitrary graph:

\begin{theorem}\label{Kimura14}{\rm\cite[Thm. 1.1]{Kimura14}} Let $G$ be a finite simple
graph. If there is a strongly disjoint family $\{B_1,\ldots,B_r\}$ of complete
bipartite subgraphs of $G$, then $\beta_{|\sigma|-r,\sigma}(R/I(G)) \neq 0$,
where $\sigma = V(B_1)\cup\cdots\cup V(B_r)$.\end{theorem}

A graph is said to be a \emph{chordal\/} if its only possible induced cycles are
triangles. In \cite{Kimura11}, Kimura proved that condition in
Theorem~\ref{Kimura14} is also a necessary condition for chordal graphs. In this
paper we show that Kimura's condition is also a necessary one for the larger
class of weakly chordal graphs (Theorem~\ref{mainwchg} below.) For general
theoretic notions we follow~\cite{Harary} and \cite{HerzogH}.

\section{Main result}

For a graph $G$, its \emph{complement graph}, denoted $\overline{G}$, is the one
with the same vertex set as $G$, an a pair of vertices form an edge of
$\overline{G}$ if and only if they do not form an edge of $G$. A graph $G$ is
said to be {\it weakly chordal\,} if neither $G$ nor its complement
$\overline{G}$ contains an induced cycle with five or more vertices. Clearly a
graph is weakly chordal if and only if its complement graph is weakly chordal.
For example, chordal graphs are weakly chordal. In fact, if $\overline{G}$
contains an induced cycle of length five, then $G$ itself contains an induced
cycle of length five; and if $\overline{G}$ contains an induced cycle of length
six or more, this cycle contains an induced subgraph consisting of the disjoint
union of two edges, which implies that $G$ contains a square as induced
subgraph. A similar argument shows that a graph whose only possible induced
cycles are squares is weakly chordal; in this case, the existence of an induced
cycle of length six or more in $\overline{G}$ would imply the existence of an
induced triangle in $G$.

Two non-adjacent vertices $u$ and $v$ of a graph $G$ form a {\it two-pair} if
any chordless path joining them has exactly two edges. An edge $e=uv$ of $G$ is
called a {\it co-pair edge} if $u$ and $v$ form a two-pair in $\overline{G}$.
The graph $G$ is said to be \emph{complete} if any pair of its vertices form an
edge of $G$.

\begin{lemma}\label{twopair}{\rm(\cite{Hay})} A graph is weakly chordal if and
only if each induced subgraph either is complete or contains a two-pair of the
subgraph.\end{lemma}

We recall that for a vertex $u$ in a graph $G$, $N(u)$ denotes the set of
\emph{neighbors} of $u$, i.e $N(u):=\{v:uv\ \text{is an edge of}\ G\}$. In
addition, we write $N[u]:=N(u)\cup\{u\}$.

\medskip The next two results are key lemmas for our main result.

\begin{lemma}\label{NV}{\rm(\cite[Lemma 7.8]{NV})} Let $e=uv$ be a co-two-pair
edge of a weakly chordal graph $G$. Then there is a complete bipartite subgraph
of $G$ with vertex set $N(u)\cup N(v)$ that contains $e$ as an edge (this union
need not be the bipartite partition of the subgraph.\end{lemma}

\proof (Sketch) Let us first partition the set $N:=N(u)\cup N(v)$ in four
disjoint sets: $T=\{u,v\}$, $U=N\backslash
 N[v]=\{u_1,u_2,\ldots\}$, $V=N\backslash N[u]=\{v_1,v_2,\ldots\}$ and
 $W=N(u)\cap N(v)$. Let us partition $W$ in three
 disjoint sets: $W_{UV}$ consisting of those $w$ such that $w$ is adjacent to every
 vertex in $U\cup V$, $W_U$ consisting of those $w$ for which there is a
 $v'\in V$ such that $w$ is not adjacent to $v'$, $W_V$ consisting of those $w$
 for which there is a $u'\in U$ such that $w$ is not adjacent to $u'$.

\vspace{-.4cm}
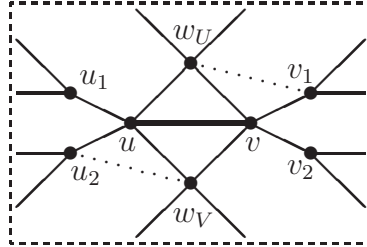
\begin{center}
\setlength{\unitlength}{.04cm} \thicklines
\begin{picture}(50,90)(20,-30)
\put(0,0){\circle*{4}} \put(0,20){\circle*{4}} \put(20,10){\circle*{4}}
\put(40,30){\circle*{4}}\put(40,-10){\circle*{4}} \put(60,10){\circle*{4}}
\put(80,0){\circle*{4}}\put(80,20){\circle*{4}}

\put(0,0){\line(2,1){20}} \put(0,20){\line(2,-1){20}}
\multiput(20,10)(.1,0){400}{\circle*{1.5}} \put(60,10){\line(2,1){20}}
\put(60,10){\line(2,-1){20}}\put(20,10){\line(1,1){38}}
\put(60,10){\line(-1,1){38}} \put(0,20){\line(-1,1){18}}
\put(0,20){\line(-1,0){18}} \put(0,0){\line(-1,0){18}}
\put(0,0){\line(-1,-1){18}} \put(80,0){\line(1,0){18}}
\put(80,0){\line(1,-1){18}}\put(80,20){\line(1,1){18}}
\put(80,20){\line(1,0){18}} \put(40,-10){\line(1,1){20}}
\put(40,-10){\line(-1,-1){18}}\put(40,-10){\line(-1,1){20}}
\put(40,-10){\line(1,-1){18}}\multiput(40,-10)(-4,1){10}{\circle*{1}}
\multiput(80,20)(-4,1){10}{\circle*{1}}

\put(16,1){$u$}\put(58,1){$v$}
\put(0,-9){$u_2$}\put(3,23){$u_1$}\put(34,38){$w_U$}\put(34,-23){$w_V$}
\put(72,25){$v_1$}\put(72,-7){$v_2$}

\put(-20,-30){\dashbox{2}(120,80)}
\end{picture}
\end{center}

\begin{itemize}
\item[1.] The vertices $u_i$ and $v_j$ are adjacent; in other case
$\overline{G}$ would contain the chordless path $u-v_j-u_i-v$, which is not
possible because $uv$ is a co-par edge of $G$.

\item[2.] $w_U\in W_U$ and $u_i$ are adjacent; in other case $\overline{G}$
would contain the chordless path $u-v'-w_U-u_i-v$, where we take $v'\in V$ such
that $w_Uv'$ is not an edge of $G$.

\item[3.] $w_V\in W_V$ and $v_j$ are adjacent.

\item[4.] If $W_U=\emptyset$, we have the bipartite partition (of $N$): $\{u\}\cup V$
and $\{v\}\cup U\cup W_V\cup W_{UV}$.

\item[5.] If $W_V=\emptyset$, we have the bipartite partition: $\{v\}\cup U$ and
$\{u\}\cup V\cup W_U\cup W_{UV}$.

\item[6.] If $W_U\neq\emptyset\neq W_V$, each $w_U\in W_U$ is adjacent to every
$w_V\in W_V$; in other case $\overline{G}$ would contain the chordless path
$u-v_j-w_V-w_U-u_i-v$.

\item[7.] If $W_U\neq\emptyset\neq W_V$ and $W_U\cup W_{UV}$ form a bipartite
subgraph, then we have the bipartite partition $\{u\}\cup V\cup W_U$ and
$\{v\}\cup U\cup W_V\cup W_{UV}$.

\item[8.] If $W_U\neq\emptyset\neq W_V$ and $W_U\cup W_{UV}$ does not form a
bipartite subgraph, then $W_V\cup W_{UV}$ form a bipartite subgraph. Hence we
have the bipartite partition $\{u\}\cup V\cup W_U\cup W_{UV}$ and $\{v\}\cup
U\cup W_V$. \qed
\end{itemize}

Recall that $\D(G)$ denotes the independence complex of the graph $G$. To
simplify notation we will write $\D G$ instead of $\D(G)$.

\begin{lemma}\label{lmesh} {\rm\cite[Claim 3.1]{Mesh}}
Let $G$ be a graph and $e=uv$ an edge of $G$. Let $G_1$ be the graph obtained
from $G$ by deleting $e$ (but not its vertices) and $G_2$ the graph obtained
from $G$ by deleting all vertices adjacent to $u$ or $v$. Then there is a long
exact sequence:
\[\cdots\ra\wt{H}_{i-1}(\Delta G_2)\ra\wt{H}_{i}(\Delta G)\ra\wt{H}_{i}(\Delta
G_1)\ra\wt{H}_{i-2}(\Delta G_2)\ra\cdots.\]
\end{lemma}

We are now ready to prove our main result.

\begin{theorem}.\label{mainwchg} Let $G$ be a weakly chordal graph with at least
one edge. If $\beta_{|\sigma|-r,\sigma}(R/I(G)) \neq 0$ for some
$\sigma\neq\emptyset$, then there exists a strongly disjoint family
$\{B_1,\ldots,B_r\}$ of complete bipartite subgraphs of $G$ for which
$\sigma=V(B_1)\cup\cdots\cup V(B_r)$.\end{theorem}

\proof We proceed by induction on the number of edges. If $G$ has just one edge,
then $\beta_{1,\{1,2\}}(R/I(G))=1$ is the only nonzero Betti number of $G$; in
this case $r=1$ and we take $B_1$ to be the edge of $G$. So we may assume that
$G$ has more than one edge.

By Hochster's formula,
\[\beta_{|\sigma|-r,\sigma}(R/I(G))=\beta_{|\sigma|-r,\sigma}(R/I_{\D G}) =
\dim\wt{H}_{r-1}((\D G)_\sigma),\] where $(\D G)_\sigma$ denotes the restriction
of $\D G$ to $\sigma$. Let $G_\sigma$ be the induced subgraph of $G$ over
$\sigma$. The graph $G_\sigma$ is weakly chordal because deleting a vertex from
a weakly chordal graph results in a graph of the same type. Since
$\D(G_\sigma)=(\D G)_\sigma$, we have that
\[\dim\wt{H}_{r-1}((\D
G)_\sigma)=\dim\wt{H}_{r-1}(\D(G_\sigma)).\] Thus, by hypothesis,
\[\beta_{|\sigma|-r,\sigma}(R/I(G_\sigma))=\dim\wt{H}_{r-1}(\D(G_\sigma)_\sigma)=
\dim\wt{H}_{r-1}(\D(G_\sigma))\neq 0.\] We will apply the last two lemmas to the
graph $G_\sigma$, so, from here on we assume that $\sigma=V(G)$ and
$G=G_\sigma$.

If $G$ is a complete graph, then $\dim\wt{H}_{r-1}(\D G)\neq 0$ implies that
$r=1$. Thus fixing a vertex of $G$ and taking as $B_1$ all edges incident with
that vertex, we are done. So, we may further assume that $G$ is not a complete
graph.

Let $e$ be a co-par edge of $G$, whose existence is ensured by
Lemma~\ref{twopair}. From Lemma~\ref{lmesh} we have the exact sequence
\[\cdots\ra\wt{H}_{r-2}(\Delta G_2)\ra\wt{H}_{r-1}(\Delta
G)\ra\wt{H}_{r-1}(\Delta G_1)\ra\wt{H}_{r-3}(\Delta G_2)\ra\cdots.\]

Since by hypothesis $\wt{H}_{r-1}(\D G)\neq 0$, in this exact sequence we have
two cases:

\medskip{\it Case} $\wt{H}_{r-1}(\Delta G_1)\neq 0$: By induction hypothesis
there is a strongly disjoint family $\{B_1,\ldots,B_r\}$ of complete bipartite
subgraphs of $G_1$. If this is not a strongly disjoint family of complete
bipartite subgraphs of $G$, necessarily there are two $3$-disjoint edges $ab\in
E(B_i)$ and $pq\in E(B_j)$ such that $e=ap$. Hence $a-q-b-p$ is a chordless path
of length three in $\overline{G}$, which contradicts that $e$ is a co-par edge
of $G$.

\medskip{\it Case} $\wt{H}_{r-2}(\Delta G_2)\neq 0$: For $\sigma_2=V(G_2)$
we have that
\[\beta_{|\sigma_2|-(r-1),\sigma_2}(R/I(G_2))=\dim \wt{H}_{r-2}(\D G_2)\neq 0.\]
By induction hypothesis there is a strongly disjoint family
$\{B_1,\ldots,B_{r-1}\}$ of complete bipartite subgraphs of $G_2$, with
$\sigma_2=V(B_1)\cup\cdots\cup V(B_{r-1})$. Because of Lemma~\ref{NV}, the
restriction of $G$ to $N(u)\cup N(v)$ contain a complete bipartite subgraph,
which we denote by $B_r$. It is clear that $\{B_1,\ldots,B_r\}$ is a strongly
disjoint family of complete bipartite subgraphs of $G$. Therefore, the proof of
the Theorem follows.\qed

\bigskip In the study of the Betti numbers $\beta_{i\sigma}(M)$, a particular
problem is the description of some homological invariants associated to them.
Two of such invariants are {\it projective dimension} and {\it regularity}:
\[\text{\rm pdim}\, M=\max\{i:\beta_{i\sigma}(M)\neq 0\};\qquad\text{\rm reg}\, M=
\max\{|\sigma|-i:\beta_{i\sigma}(M)\neq 0\}.\]

\begin{remark}\label{DalKumm} By~\cite[Prop. 5.3]{DalKumm}, the Betti numbers
of a monomial ideal over the field of rational numbers can be obtained from the
Betti numbers over any field by a sequence of consecutive cancellations. Since
the information of the Betti numbers given in Theorem~\ref{mainwchg} is entirely
combinatorial, it follows that the projective dimension and regularity of weakly
chordal graphs are independent of the characteristic of the base
field.\end{remark}

Let $\text{\rm imn}(G)$ denote the maximum number of parwise $3$-disjoint edges
of $G$. Further, following~\cite{NV}, we define $d(G) =
\max\{\sum_{i=1}^r|V(B_i)|-r\}$, where the maximum is taken over all the
strongly disjoint families $\{B_1,\ldots,B_r\}$ of complete bipartite subgraphs
of $G$.

\begin{corollary}\label{maincor} For any weakly chordal graph $G$ it holds that
\begin{itemize}
 \item[$(i)$] {\rm\cite[Prop. 20]{Wood}}\,\ \ $\text{\rm reg}\, R/I(G)=\text{\rm imn}(G)$.
 \item[$(ii)$] {\rm\cite[Thm. 7.7]{NV}}\ \ $\text{\rm pdim}\, R/I(G)=d(G)$.
\end{itemize}\end{corollary}

\proof (i) Let $i$ and $\sigma$ be such that $\beta_{i\sigma}(R/I(G))\neq 0$ and
$\text{\rm reg}\, R/I(G)=|\sigma|-i=:r$. Hence
$\beta_{|\sigma|-r,\sigma}(R/I(G))\neq 0$. By Theorem~\ref{mainwchg} there
exists a strongly disjoint family $\{B_1,\ldots,B_r\}$ of complete bipartite
subgraphs of $G$. Thus $r\leq\text{\rm imn}(G)$. The other inequality,
$\text{\rm imn}(G)\leq r$, is well-known \cite{Katzman}.

(ii) Let $\{B_1,\ldots,B_r\}$ be a strongly disjoint family of complete
bipartite subgraphs of $G$ such that $m=\sum_{i=1}^r|V(B_i)|-r$ is maximum. Set
$\sigma=V(B_1)\cup\cdots\cup V(B_r)$ and $p:=\text{\rm pdim}\, R/I(G)$. From
Theorem~\ref{Kimura14} we have that $\beta_{|\sigma|-r,\sigma}(R/I(G))\neq 0$,
hence $p\geq m$. On the other hand, by definition of $p$, there is a $\tau$ such
that $\beta_{p\tau}(R/I(G))\neq 0$. Write $p=|\tau|-q$ for some $q$. By
Theorem~\ref{mainwchg} there is a strongly disjoint family
$\{B_1',\ldots,B_q'\}$ of complete bipartite subgraphs of $G$ such that
$\tau=V(B_1')\cup\cdots\cup V(B_q')$. Then, by maximality of $m$, it follows
that $p=|\tau|-q\leq m$. \qed

\begin{remark} It is well-known that the projective dimension of $R/I(G)$ equals the
big-height of $I(G)$ when $G$ is a chordal graph~\cite[Cor. 5.6]{Kimura11},
\cite[Cor. 3.33]{MorVila}, \cite[Thm. 3.2]{ChFVT}. This result cannot be
extended to weakly chordal graphs. For instance, a square has projective
dimension $3$, while its big-height equals $2$.\end{remark}

\begin{remark} There are graphs containing pentagons as induced subgraphs and its projective
dimension and regularity depend on the characteristic of the base field. For
instance, for Katzman's example~\cite{Katzman}
 
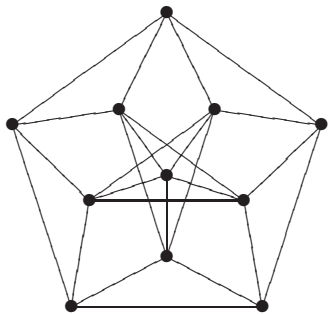
\vspace{-.2cm}
 \[
 \xygraph{
 !{<0cm,0cm>:<0cm,1.2cm>:<-1.2cm,0cm>}
 !{(0,0);a(0)**{}?(0)}*{\bullet}@\cir{}="a0"
 !{(0,0);a(-36)**{}?(.9)}*{\bullet}@\cir{}="a1"
 !{(0,0);a(36)**{}?(.9)}*{\bullet}@\cir{}="a2"
 !{(0,0);a(108)**{}?(.9)}*{\bullet}@\cir{}="a3"
 !{(0,0);a(180)**{}?(.9)}*{\bullet}@\cir{}="a4"
 !{(0,0);a(252)**{}?(.9)}*{\bullet}@\cir{}="a5"
 !{(0,0);a(0)**{}?(1.8)}*{\bullet}@\cir{}="b1"
 !{(0,0);a(72)**{}?(1.8)}*{\bullet}@\cir{}="b2"
 !{(0,0);a(144)**{}?(1.8)}*{\bullet}@\cir{}="b3"
 !{(0,0);a(216)**{}?(1.8)}*{\bullet}@\cir{}="b4"
 !{(0,0);a(288)**{}?(1.8)}*{\bullet}@\cir{}="b5"
 "a0"-"a1" "a0"-"a2" "a0"-"a3" "a0"-"a4" "a0"-"a5"
 "a1"-"a3" "a3"-"a5" "a5"-"a2" "a2"-"a4" "a4"-"a1"
 "b1"-"b2" "b2"-"b3" "b3"-"b4" "b4"-"b5" "b5"-"b1"
 "a1"-"b1" "a2"-"b2" "a3"-"b3" "a4"-"b4" "a5"-"b5"
 "a1"-"b5" "a5"-"b4" "a4"-"b3" "a3"-"b2" "a2"-"b1"
 }\]

\noindent these invariants are $8$ and $2$, respectively in characteristic $\neq
2$, while they are $9$ and $3$, respectively in characteristic $2$.
\end{remark}


\end{document}